\documentclass[12pt]{article}

\usepackage{verbatim}
\usepackage{amssymb}  
\usepackage{amsmath}
\usepackage{amsthm}
\usepackage{mathrsfs}
\usepackage[colorlinks]{hyperref}
\usepackage{url}
\usepackage{fullpage}
\usepackage{tikz}

\usepackage{authblk}
\usepackage{graphicx}
\usetikzlibrary{backgrounds,fit,decorations.pathreplacing}
\usetikzlibrary{arrows,shapes,positioning}
\usetikzlibrary{calc,through,backgrounds,chains}
\usetikzlibrary{arrows,shapes,snakes,automata,backgrounds,petri}

\newtheorem{Theorem}[equation]{Theorem}

\newtheorem{Corollary}[equation]{Corollary}
\newtheorem{Lemma}[equation]{Lemma}
\newtheorem{Proposition}[equation]{Proposition}

\theoremstyle{definition}

\newtheorem{Example}[equation]{Example}

\theoremstyle{remark}

\newtheorem{Remark}[equation]{Remark}

\numberwithin{equation}{section}
\numberwithin{figure}{section}

\newcommand{\C}{\mathbb{C}}

\newcommand{\Z}{\mathbb{Z}}

\begin{document}

\title{The Rook Monoid is Lexicographically Shellable}
\date{\today}
\author[1]{Mahir Bilen Can}
\affil[1]{{\small Tulane University, New Orleans; mahirbilencan@gmail.com}}

\maketitle

\begin{abstract}
We prove that the Bruhat-Chevalley-Renner order on the rook monoid 
is EL-shellable. We determine the homeomorphism type of the associated order complex. 
\vspace{.2cm}

\noindent 
\textbf{Keywords:}  Rook monoid, Bruhat-Chevalley-Renner order, EL-shellability, homeomorphism type.\\ 
\noindent 
\textbf{MSC: 05E45, 20M32} 
\end{abstract}

\section{\textbf{Introduction}}

Let $n$ be a positive integer, and let $[n]$ denote the set $\{1,\dots,n\}$.
Let $P$ be a finite graded poset of rank $n$, with minimum 
and maximum elements denoted by $\hat{0}$ and $\hat{1}$, respectively.
We denote by $C(P)$ the set of pairs $(x,y)$ from $P \times P$ 
such that $y$ covers $x$. 
The poset $P$ is called {\em lexicographically shellable}, or {\em EL-shellable}, 
if there exists a map $f: C(P) \rightarrow [n]$ such that
\begin{enumerate}
\item[(1)] in every interval $[x,y]\subseteq P$, there exists 
a unique maximal chain 
$\mathfrak{c}:\ x=x_0 < x_1 < \cdots < x_{k+1} =y$ 
such that $f(x_i,x_{i+1}) \leq f(x_{i+1},x_{i+2})$ for $i=0,\dots,k-1$;
\item[(2)] the sequence $f(\mathfrak{c}):=(f(x,x_1),\dots,f(x_k,y))$ 
of the unique chain $\mathfrak{c}$ of (1) is lexicographically 
first among all sequences of the form $f(\mathfrak{d})$, 
where $\mathfrak{d}$ is a maximal chain in $[x,y]$.
\end{enumerate}
Introduced by Bj\"{o}rner in \cite{Bjorner80}, 
the notion of lexicographic shellability has important 
topological consequences. For example, it is known that if 
the poset is lexicographically shellable, then the
order complex of the poset is a shellable simplicial complex.

In this paper, we are concerned with the question 
of lexicographic shellability of the Bruhat-Chevalley-Renner order on the monoid
of injective partial transformations, which is called the {\em rook monoid}, and 
denoted by $R_n$. By definition, $R_n$ is the finite monoid consisting of 
$n\times n$ $0/1$ matrices with at most one 1 in each row and in each column. 
To define the Bruhat-Chevalley-Renner order, we first explain the role of matrices for $R_n$.

Let $M_n$ denote the algebraic monoid of all $n\times n$ matrices with entries in $\C$,
and let $GL_n$ denote its unit-group.
We denote by $B_n$ the Borel subgroup, which consists of  
invertible upper triangular elements from $GL_n$. We consider 
the action 
\begin{align}
(B_n\times B_n) \times M_n &\longrightarrow M_n \notag \\
((b_1,b_2),m) &\longmapsto b_1m b_2^{-1} \label{I:A:action}
\end{align} 
By Bruhat-Chevalley decomposition, the orbits of the restriction of this action to $GL_n$
are parametrized by the symmetric group on $n$ letters, denoted by $S_n$. 
Then the {\em Bruhat-Chevalley} order on $S_n$ is defined by 
\begin{align}\label{A:BC on S_n}
x \leq y \iff B_n x B_n \subseteq \overline{B_n y B_n} ,\ \text{ where $x,y\in S_n$.}
\end{align}
Here, the bar on the orbit indicates the Zariski closure in $GL_n$.  
More generally, the decomposition of $M_n$ into $B_n\times B_n$-orbits is parametrized
by the rook monoid,  
$$
M_n = \bigsqcup_{r\in R_n} B_n r B_n,
$$
and the {\em Bruhat-Chevalley-Renner order} on $R_n$ is defined by 
\begin{align}\label{A:BC on R_n}
r \leq t \iff B_n r B_n \subseteq \overline{B_n t B_n} ,\ \text{ where $r,t\in R_n$,}
\end{align}
see~\cite{Renner86}.
In (\ref{A:BC on R_n}), the bar indicates the Zariski closure in $M_n$.

We are now ready to state our main results, and 
explain how we structured our paper. 
In Section~\ref{S:definitions}, we introduce our 
notation and review in more detail the Bruhat-Chevalley-Renner  
order on the rook monoid as well as on the symmetric group.
In particular, we review a result of Edelman, which states that 
the Bruhat-Chevalley order on symmetric groups is lexicographically 
shellable. 
In Section~\ref{S:labelingoftherooks},
by describing an explicit labeling of the edges of the Hasse diagram
of the Bruhat-Chevalley-Renner order on $R_n$, we state 
our main result. Its proof will be given in Section~\ref{S:proofs}.
\begin{Theorem}
If $n$ is a positive integer, then $(R_n,\leq)$ is an EL-shellable poset.
\end{Theorem}

The {\em order complex} of a poset $P$ is 
the abstract simplicial complex whose closed
sets are the finite chains in $P$. 
In the final section of our paper, that is Section~\ref{S:Final}, 
we consider the order complexes of intervals in $R_n$. 
We finish our paper by proving the second main result of our paper.
\begin{Theorem}
The order complex of an open interval in
$R_n$ is homeomorphic to a sphere or to a ball. 
In particular, $\Delta(R_n)$ triangulates a ball of dimension $n^2$.
\end{Theorem}

\noindent \textbf{Acknowledgements.}
We thank Lex Renner for teaching us the theory of linear algebraic monoids. 
We thank the referees whose comments and suggestions greatly improved the quality of our paper. 
In particular, we are grateful to the referee who showed us another EL-labeling on the rook monoid. 
Finally, we thank Michelle Wachs and Warut Thawinrak. 
\section{\textbf{Background}}\label{S:definitions}

\subsection{Reductive monoids and Bruhat-Chevalley-Renner order.}\label{SS:reductivemonoids}

Let $G$ be a connected reductive complex algebraic group. 
We fix a maximal torus $T$, and a Borel subgroup $B$ in $G$ 
such that $T\subset B\subset G$. Then the {\em Weyl group} of $(G,T)$ 
is the quotient group $W:=N_G(T)/T$, where $N_G(T)$ 
is the normalizer of $T$ in $G$. Note that for the pair $(GL_n,T)$, where $T$ is the maximal torus 
of diagonal matrices in $GL_n$, the Weyl group is isomorphic to 
the symmetric group, $S_n$.
The Bruhat-Chevalley order on a Weyl group $W$ is defined by  
$w \leq v \iff Bw B \subseteq \overline{B v B}\ \text{  for } v,w\in W$.
It is well-known that such orders are lexicographically shellable, see \cite{BW82, Edelman81,Proctor82}.

An algebraic monoid is an algebraic variety $M$ 
together with an associative multiplication morphism $m: M\times M \rightarrow M$ such that
 there is a neutral element.
For an algebraic monoid $M$, the {\em unit-group}, that is the group of invertible elements of $M$,
has the structure of an algebraic group.  
A {\em reductive monoid} is an algebraic monoid $M$ such that the unit-group of $M$ 
is a connected reductive algebraic group.

Interesting examples of reductive monoids are built from representations.
Let $\rho: G_0 \rightarrow GL(V)$ be a rational representation of a semisimple algebraic group $G_0$. 
Let us denote by $\C^*$ the multiplicative group of nonzero scalar matrices in $End(V)$. 
Then the Zariski closure $M:=\overline{\C^*\cdot \rho(G_0)}$ in $End(V)$ 
has the structure of a reductive monoid. The normality of such monoids is intimately related 
to the highest weight of the representation, see~\cite{DeConcini}.

Let $G$ denote the unit-group of a reductive monoid $M$, 
and let $T \subset B$ be a maximal 
torus and a Borel subgroup in $G$. 
It is shown in \cite{Renner86} that there is a decomposition of $M$ 
into double cosets of $B$, $M = \bigsqcup_{r\in R} B \dot{r} B$,
where $R:=\overline{N_G(T)}/T$ and $\dot{r}\in \overline{N_G(T)}$.
The parametrizing object $R$, which is a finite monoid, 
is called the {\em Renner monoid} of $M$. 
The Bruhat-Chevalley-Renner order on $R$ is defined by 
$$
w \leq v \iff Bw B \subseteq \overline{B v B}\ \text{  for } v,w\in R.
$$

If $\rho: G_0 \rightarrow GL(K^n)$ is the defining representation of $G_0=SL_n$, 
then the Renner monoid $R$ is isomorphic to the rook monoid $R_n$. 
The Weyl group $W$ of $(G,T)$ is the unit-group of $R$, 
and the Bruhat-Chevalley order on $W$ extends to the Bruhat-Chevalley-Renner order on $R$.

There is a cross section lattice of idempotents, denoted by $\Lambda$, in $R$ such that 
\begin{equation*}
M= \bigsqcup_{e\in \Lambda} GeG \ \text{ and }\ R= \bigsqcup_{e\in \Lambda} WeW.
\end{equation*}
It is shown by Putcha in \cite{Putcha02} that the diagonal subposet 
$WeW\subseteq R$ ($e\in \Lambda$) is a retract of a shellable poset.
Note that, if $\Lambda\setminus \{0\}$ has a unique minimal element, for example   
when $M$ is the Zariski closure of an irreducible representation of a connected reductive group, 
then the cross section lattice is an (upper) semimodular lattice. In this case, $\Lambda$ is shellable. 
We expect that the Renner monoid of a normal reductive monoid is EL-shellable.

For more information about linear algebraic monoids, we recommend the resources 
~\cite{Putcha88,Renner05,Solomon95}.

\subsection{Lexicographic shellability}\label{SS:lexicographic}

Let $P$ be a graded poset of rank $n$ with 
a maximum and a minimum element, denoted by $\hat{1}$ and $\hat{0}$ respectively. 
Then all maximal chains of $P$ have equal length $n$. 
If $x$ and $y$ are two elements from $P$ such that $y$ covers $x$, 
then we will write $x \lessdot y$. 
We denote by $C(P)$ the set of covering relations in $P$,  
$C(P)= \{(x,y)\in P\times P:\ x \lessdot y \}$.
An \textit{edge-labeling} on $P$ is a map of the form 
$f=f_{P,\varGamma}: C(P) \rightarrow \varGamma$, where $\varGamma$ is 
a totally ordered set.
The  \textit{Jordan-H\"{o}lder sequence} (with respect to $f$) 
of a maximal chain $\mathfrak{c}: x_0 \lessdot x_1\lessdot \cdots \lessdot x_{n-1}\lessdot x_n$ of 
$P$ is the $n$-tuple 
$f(\mathfrak{c})= (f(x_0,x_1), f(x_1,x_2),\dots, f(x_{n-1},x_n)) \in \varGamma^n$.
We fix an edge labeling, denoted $f$, 
and we fix a maximal chain $\mathfrak{c}:\ x_0< x_1 < \cdots < x_n$.
We call both of the maximal chain $\mathfrak{c}$ and its image $f(\mathfrak{c})$ {\em increasing} 
if the following inequalities hold true:
$$
f(x_0,x_1) \leq f(x_1,x_2) \leq \cdots \leq f(x_{n-1},x_n).
$$

Let $k$ be a positive integer. 
We consider the lexicographic order on the $k$-fold 
cartesian product, 
$\varGamma^k = \varGamma \times \cdots \times \varGamma $. 
An edge labeling $f:C(P) \rightarrow \varGamma$ 
is called an {\em EL-labeling} if
\begin{enumerate}
\item in every interval $[x,y] \subseteq P$ of rank $k>0$ there exists a unique maximal chain $\mathfrak{c}$ such that 
$f(\mathfrak{c}) \in \varGamma^k$ is increasing,
\item the Jordan-H\"{o}lder sequence $f(\mathfrak{c}) \in \varGamma^k$ of the unique chain $\mathfrak{c}$ from (1) is 
the smallest among the Jordan-H\"{o}lder sequences of maximal chains $x=x_0\lessdot x_1 \lessdot \cdots \lessdot x_k=y$.
\end{enumerate}
A poset $P$ is called {\em EL-shellable} if it has an EL-labeling.

\begin{Remark}
There are various lexicographic shellability conditions in the literature, 
and the EL-shellability that is defined here is among the stronger ones, see~\cite[Section 3.2]{Wachs07}.
A deep relationship between EL-shellability of a Coxeter group $W$ and the Kazhdan-Lusztig theory of the Hecke algebra 
associated with $W$ is found by Dyer in \cite{Dyer93}.
\end{Remark}

\subsection{The symmetric group}\label{SS:symmetricgroup}
We already introduced the notation $S_n$ for the set of all permutations of $[n]$. 
We represent the elements of $S_n$ in one-line notation, $w=(w_1,\dots,w_n)\in S_n$, 
so that $w(i)=w_i$. 
It is well-known that $S_n$ is a graded poset with respect to Bruhat-Chevalley order. 
Let $B_n$ denote, as before, the subgroup of invertible upper triangular matrices in $GL_n$. 
The grading on $S_n$ is given by the {\em length function}
\begin{equation}\label{E:lengthsymmetric}
\ell(w)= \dim (B_n w B_n) -\dim(B_n)=inv(w),
\end{equation}
where
\begin{equation}\label{E:invpermutation}
inv(w) = |\{ (i,j):\ 1\leq i < j \leq n,\ w_i>w_j \}|.
\end{equation}
Note that $\dim B_n = {n+1\choose 2}.$

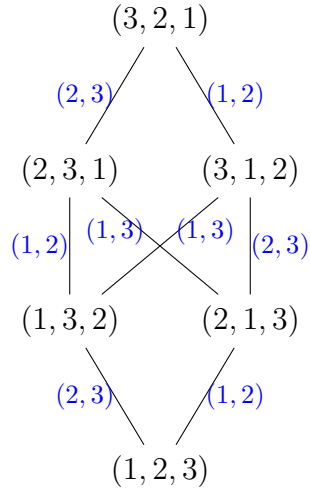
\begin{figure}[htp]
\begin{center}
\scalebox{.85}{
\begin{tikzpicture}[scale=.4]
media/.style={font={\footnotesize}},

\node at (0,0) (a) {$(1,2,3)$};
\node at (-3,5) (b1) {$(1,3,2)$};
\node at (3,5) (b2) {$(2,1,3)$};
\node at (-3,10) (c1) {$(2,3,1)$};
\node at (3,10) (c2) {$(3,1,2)$};
\node at (0,15) (d) {$(3,2,1)$};
\draw[-] (a) -- (b1);
\draw[-] (a) -- (b2);
\draw[-] (b1) -- (c1);
\draw[-] (b1) -- (c2);
\draw[-] (b2) -- (c1);
\draw[-] (b2) -- (c2);
\draw[-] (c1) -- (d);
\draw[-] (c2) -- (d);
\begin{scope}[blue]
\node at (-2.5,2.5) {\footnotesize{$(2,3)$}};
\node at (2.5,2.5) {\footnotesize{$(1,2)$}};
\node at (-4,7.5) {\footnotesize{$(1,2)$}};
\node at (4,7.5) {\footnotesize{$(2,3)$}};
\node at (-1.5,8) {\footnotesize{$(1,3)$}};
\node at (1.5,8) {\footnotesize{$(1,3)$}};
\node at (-2.5,12.5) {\footnotesize{$(2,3)$}};
\node at (2.5,12.5) {\footnotesize{$(1,2)$}};
\end{scope}

\end{tikzpicture}
}
\end{center}
\caption{EL-labeling  of $S_3$}
\label{F:diagram}
\end{figure}

The Bruhat-Chevalley order on $S_n$ is the partial order 
generated by the transitive closures of the following covering relations: 
$x=(a_1,\dots.,a_n)$ is covered by $y=(b_1,\dots,b_n)$ 
if  $\ell(y)=\ell(x)+1$, and
\begin{enumerate}
\item $a_k=b_k$ for $k\in \{ 1,\dots,\widehat{i},\dots,\widehat{j},\dots,n\}$ (hat means omit those numbers), 
\item $a_i=b_j$, $a_j=b_i$, and $a_i<a_j$.
\end{enumerate}
Let $\varGamma=[n] \times [n]$ denote the poset of pairs, ordered lexicographically: 
$(i,j) \leq (r,s)$ if $i<r$, or  $i=r$ and $j<s$. 
For permutations $x,y\in S_n$ given as above, define $f(x,y):=(a_i,a_j)$.
In Figure \ref{F:diagram}, we depict the corresponding EL-labeling for $S_3$.

\begin{Theorem} (\cite{Edelman81})
Let $f$ denote the edge-labeling on the Bruhat-Chevalley order on $S_n$ as defined in the previous paragraph.  
Then $S_n$ is an EL-shellable poset with respect to $f$. 
\end{Theorem}

\subsection{The rook monoid}\label{SS:rookmonoid}

The rook monoid $R_n$ is a graded poset (see \cite{Renner86}) and its rank
function is given by 
\begin{equation*}
\ell( x )= \dim (B_n x B_n),\  x \in R_n.
\end{equation*}
There is a combinatorial formula for $\ell(x)$, $x\in R_n$ that is similar to (\ref{E:lengthsymmetric}).
To explain, we represent the elements of $R_n$ by $n$-tuples.
Let $x=(x_{ij}) \in R_n$, and we define the sequence $(a_1,\dots,a_n)$ by
\begin{equation}\label{E:oneline}
a_j = 
\begin{cases}
0  &\text{if the $j$-th column consists of zeros,}\\
i &\text{if $x_{ij}=1$.}
\end{cases}
\end{equation}
By abuse of notation, we denote both of the matrix and the corresponding 
sequence $(a_1,\dots,a_n)$ by $x$.
For example, the associated sequence of the partial permutation matrix 
\begin{equation*}
x=\begin{pmatrix}  
0 & 0 & 0 & 0 \\
0 & 0 & 0 & 0 \\
1 & 0 & 0 & 0 \\ 
0 & 0 & 1 & 0
\end{pmatrix}
\end{equation*}
is $x=(3,0,4,0)$.

Let $x=(a_1,\dots.,a_n)$ be an element from $R_n$.  
A pair of indices $(i,j)$ with $1\leq i<j \leq n$ is called a {\em coinversion pair} for $x$ 
if $0< a_i < a_j$. We denote the number of coinversion pairs of $x$ by $coinv(x)$. 
\begin{Example}
If $x$ is given by $x=(4,0,2,3)$, then the only coinversion pair for $x$ is $(3,4)$. 
\end{Example}

In \cite{CR11}, it is shown that the dimension $\ell(x)=\dim(B_nxB_n)$ of an orbit $B_n x B_n$, $x\in R_n$ 
is given by
\begin{equation}\label{E:dimforminv}
\ell(x)  = \left(\sum_{i=1}^n a_i^*\right)  - coinv(x),\ \text{where}\
a_i^* =
\begin{cases}
a_i+n-i,  & \text{if}\  a_i\neq 0 \\
0,   & \text{if}\ a_i=0
\end{cases}
\end{equation}
We reformulate (\ref{E:dimforminv}) as follows:
\begin{Proposition}\label{P: sum inv}
Let $x=(a_1,\dots, a_n)\in R_n$. Then we have 
$\ell( x) = \sum a_i + inv (x)$,
where $inv(x) = |\{(i,j):\ 1\leq i< j \leq n,\ a_i>a_j \} |$.
\end{Proposition}
\begin{proof}
Since $coinv(x) = \sum_{i=1}^n c_i$, where $c_i$ is the number of $j \in [n]$ such that 
$i < j$ and $0< a_i < a_j$, we can rewrite (\ref{E:dimforminv}) as follows:
\begin{align*}
\ell(x) &= \sum_{i=1}^n (a_i^*- c_i) =\sum_{i=1}^n (a_i + d_i^*) =\sum_{i=1}^n a_i + \sum_{i=1}^n d_i^*, 
\end{align*}
where
\begin{equation*}
d_i^*= 
\begin{cases}
n-i - c_i,  & \text{if}\  a_i\neq 0, \\
0,   & \text{if}\ a_i=0.
\end{cases}
\end{equation*}
Observe now that $d_i^*$ is equal to the number of $j > i$ such that $a_i > a_j \geq 0$. 
This finishes the proof.
\end{proof}

The following consequence of Proposition \ref{P: sum inv} agrees with equation (\ref{E:lengthsymmetric}).
\begin{Corollary}\label{C:length for Sn}
Let $w=(a_1,\dots,a_n) \in S_n$ be a permutation.
Then the length of $w$ as an element of the poset $(R_n,\leq)$ is given by 
$\ell(w) = {n+1 \choose 2}+inv(w)$. 
\end{Corollary}

The first concrete description of the Bruhat-Chevalley-Renner order on $R_n$ is given in 
\cite[Theorem 3.8]{PPR97}. We state it here for convenience. 
\begin{Theorem}\label{T:PPR} 
Let $x = (a_1,\dots,a_n)$, $y=(b_1,\dots,b_n) \in R_n$. The Bruhat-Chevalley order on 
$R_n$ is the smallest partial order on $R_n$ generated by declaring 
$x \leq y$ if either 
\begin{enumerate}
\item there exists an $1 \leq i \leq n$ such that $b_i> a_i$ and $b_j = a_j$ for all $j\neq i$, or
\item  there exist  $1 \leq i < j \leq n$ such that $b_i=a_j,\ b_j=a_i$ with $b_i > b_j$, 
and for all $k\notin \{i,j\}$, $b_k = a_k$.
\end{enumerate}
\end{Theorem}

The covering relations of the order are analyzed in detail in \cite{CR11}, 
and the following two lemmas are useful in this analysis.
The original statement of the first lemma in~\cite{CR11} had a minor gap,
which was pointed out to us by  Michelle Wachs. 
We thank Warut Thawinrak, who helped us to state the lemma in a better way. 
For completeness, we present the corrected version of the lemma and its proof here.

\begin{Lemma} \label{L:PPRcovering0}
Let $x=(a_1,\dots,a_n)$ and $y = (b_1,\dots,b_n)$ be elements of $R_n$. 
We assume that $a_k=b_k$ for all $k =\{ 1,\dots,\widehat{i},\dots,n\}$ and that $a_i < b_i$. 
Then, $\ell(y) = \ell( x)+ 1$ if and only if either
\begin{enumerate}
\item $0=a_i$, and there exists a sequence of indices $1 \leq j_1< \cdots < j_s < i$ 
such that the set $\{a_{j_1},\dots, a_{j_s}\}$ is equal to $\{1,\dots, s\}$, $b_i = s+1$, and $b_j=a_j > 0$ 
for $j>i$, or 
\item $0< a_i$, and there exists a sequence of indices $1 \leq j_1< \cdots < j_s < i$ 
such that the set $\{a_{j_1},\dots,a_{j_s}\}$ is equal to 
$\{a_i+1,\dots,a_i+s	\}$, and $b_i=a_i+s+1$.
\end{enumerate}
\end{Lemma}

\begin{proof}
Let $x,y\in R_n$ be two rooks as in the hypothesis, and we assume that 
$\ell(y)= \ell(x) +1$. By Theorem~\ref{T:PPR}, we know that $x\lessdot y$. 
Obviously, either $b_i = a_i +1$, or $b_i > a_i +1$. In the former case, we have two possibilities:
$a_i = 0$, or $a_i > 0$. In the former case, it follows from Proposition \ref{P: sum inv} that, 
unless $a_j=b_j > 0$ for all $j>i$, the equality $\ell(y) = \ell(x) +1$ is not true. 
In the latter case, there is nothing to say. 
In conclusion, if $b_i = a_i +1$, then 1. holds true. On the other hand, 
if it is true that $b_i = a_i + d$ for some $d > 1$, then more analysis is needed;
\begin{eqnarray*}
\ell(y) &=&  \sum_{j=1}^n b_j^* - coinv(y)\\
   &=&  ( \sum_{j=1, j\neq i}^n a_j^* )+ b_i^*  - coinv(y)\\
   &=&  ( \sum_{j=1, j\neq i}^n a_j^* )+ a_i+d+n-i  - coinv(y)\\
   &=& ( \sum_{j=1}^n a_j^* )+ d - coinv(y)\\
   &=& \ell(x)+d+coinv(x)-coinv(y).
\end{eqnarray*}
Hence, $d+coinv(x)-coinv(y)=1$, or $coinv(y)-coinv(x)=d-1$.
We inspect the difference $coinv(x)-coinv(y)$.
If $(k,i)$ with $k<i$ is a coinversion for $x$, then it stays to be a coinversion for $y$.
Clearly, this is also true for the pairs of the form $(k,l)$ where $k<i<l$, or $i<k<l$, or $k<l<i$.
Therefore, the difference between $coinv(y)$ and $coinv(x)$ occurs at the pairs of the form
\begin{enumerate}
\item $(k,i),\ k<i$ such that $a_i < a_k < b_i$, or
\item $(i,l),\ i<l$, such that $a_i < a_l < b_i$.
\end{enumerate}
In the first case, some new coinversions are added, and in the second case some coinversions are deleted.
Let us call the number of pairs of the first type by $n_1$, 
and we call the number of pairs of the second type by $n_2$. 
Then, $coinv(y) = coinv(x) + n_1-n_2$, or $coinv(y)-coinv(x) = n_1-n_2$. 
We know that $0 \leq n_1, n_2 \leq d-1$ since $b_i=a_i+d$. 
Hence, we have that $n_1=d-1$, and that $n_2=0$. 
Therefore, the following is true: 
Any $a_k$ between $a_i$ and $a_i+d=b_i$ appears before the $i$-th position. 
This completes the proof of ``$\Rightarrow$'' direction of our lemma.

Next, we prove the ``$\Leftarrow$'' direction.
If 1. holds true, then $b_i = a_i +1$. In this case, it is straightforward to check that $\ell(y) = \ell(x) +1$. 
So, we assume that there exists a sequence of indices, $1 \leq j_1< \cdots < j_s < i$, 
such that the set $\{a_{j_1},\dots,a_{j_s}\}$ is equal to $\{a_i+1,\dots,a_i+s\}$, and $b_i=a_i+s+1$.  
Then
\begin{eqnarray*}
\ell(y) &=&  \sum_{j=1}^n b_j^* - coinv(y)\\
   &=&  ( \sum_{j=1, j\neq i}^n a_j^* )+ b_i^*  - coinv(y)\\
   &=&  ( \sum_{j=1, j\neq i}^n a_j^* )+ a_i+s+1+n-i  - coinv(y)\\
   &=& ( \sum_{j=1}^n a_j^* )+ s+1 - coinv(y).
\end{eqnarray*}
Now it suffices to show that $coinv(y)=s+coinv(x)$. 
Observe that, when we replace $a_i$ by $b_i$, the set of pairs, $\{ (j_k, i)|\ k=1,\dots,s\}$, 
whose elements are non-coinversion pairs for $x$, become coinversion pairs for $y$. 
Also, if we replace the entry $a_i$ by $b_i$, then the coinversion pairs for $x$ of the form $(l,i)$ or $(i,l)$ (where $l\neq j_k$) 
are coinversion pairs for $y$ as well. Therefore, we see that  
$coinv(y) = s+ coinv(x)$, hence that $\ell(y) = \ell(x) +1$. This finishes the proof of our lemma.
\end{proof}

\begin{Example}
Let $x=(4,0,5,0,3,1)$ and $y=(4,0,5,0,6,1)$. Then $\ell(x)= 21$ and $\ell(y)=22$. 
If $z=(4,0,5,0,3,2)$, then $\ell(z)=22$.
\end{Example}

\begin{Lemma} \label{L:PPRcovering1}
Let $x=(a_1,\dots,a_n)$ and $y = (b_1,\dots,b_n)$ be two elements of $R_n$. 
Suppose that $a_j=b_i,\ a_i = b_j$ and $b_j < b_i$ where $i < j$. 
Furthermore, suppose that for all $k\in \{1,\dots\widehat{i},\dots,\widehat{j},\dots,n\}$, $a_k=b_k$.
Then, $\ell(y) = \ell( x)+ 1$ if and only if for $s=i+1,\dots,j-1$, either $a_j< a_s$, or $a_s < a_i$.
\end{Lemma}

\begin{Example}
Let $x=(2,6,5,0,4,1,7)$ and $y=(4,6,5,0,2,1,7)$. 
Then $\ell(x)= 35$ and $\ell(y)=36$. If $z=(7,6,5,0,4,1,2)$, then $\ell(z)=42$.
\end{Example}

\section{\textbf{An EL-labeling  of $R_n$}}\label{S:labelingoftherooks}

Let $x$ and $y$ be two elements from $R_n$ such that $x\lessdot y$. 
We call $x \lessdot y$ a {\em type 1 covering relation} if it is as in Lemma \ref{L:PPRcovering0}, 
and we call it a {\em type 2 covering relation} if it is as in Lemma \ref{L:PPRcovering1}.
Let $\varGamma$ denote the finite set $\varGamma :=\{0,1,\dots,n\} \times \{1,\dots,n\}$,
totally ordered with respect to lexicographic ordering.
We define $F:C(R_n)\longrightarrow \varGamma$ by 
\begin{equation}\label{E:rooklabel}
F(x,y)=
\begin{cases}
(a_i,b_i),  & \text{if}\ y\ \text{covers}\ x\ \text{by type}\ 1\\
(a_i,a_j),  & \text{if}\ y\ \text{covers}\ x\ \text{by type}\ 2.
\end{cases}
\end{equation}

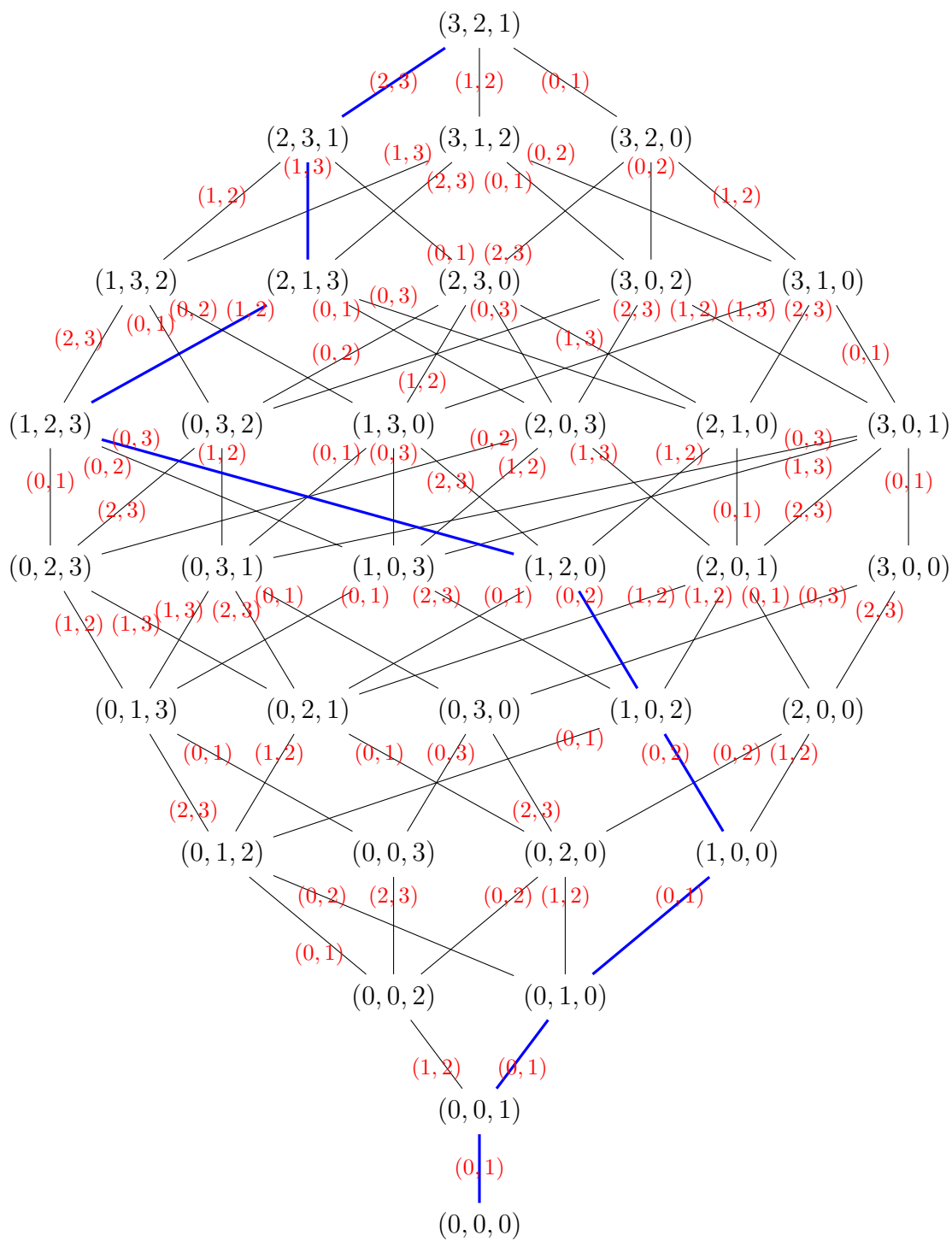
\begin{figure}[htp]
\begin{center}

\scalebox{.65}{
\begin{tikzpicture}[scale=.5]

\node at (0,2) (a) {$(0,0,0)$};

\node at (0,6) (b) {$(0,0,1)$};

\node at (-3,10) (c1) {$(0,0,2)$};
\node at (3,10) (c2) {$(0,1,0)$};

\node at (-9,15) (d1) {$(0,1,2)$};
\node at (-3,15) (d2) {$(0,0,3)$};
\node at (3,15) (d3) {$(0,2,0)$};
\node at (9,15) (d4) {$(1,0,0)$};

\node at (-12,20) (e1) {$(0,1,3)$};
\node at (-6,20) (e2) {$(0,2,1)$};
\node at (0,20) (e3) {$(0,3,0)$};
\node at (6,20) (e4) {$(1,0,2)$};
\node at (12,20) (e5) {$(2,0,0)$};

\node at (-15,25) (f1) {$(0,2,3)$};
\node at (-9,25) (f2) {$(0,3,1)$};
\node at (-3,25) (f3) {$(1,0,3)$};
\node at (3,25) (f4) {$(1,2,0)$};
\node at (9,25) (f5) {$(2,0,1)$};
\node at (15,25) (f6) {$(3,0,0)$};

\node at (-15,30) (g1) {$(1,2,3)$};
\node at (-9,30) (g2) {$(0,3,2)$};
\node at (-3,30) (g3) {$(1,3,0)$};
\node at (3,30) (g4) {$(2,0,3)$};
\node at (9,30) (g5) {$(2,1,0)$};
\node at (15,30) (g6) {$(3,0,1)$};

\node at (-12,35) (h1) {$(1,3,2)$};
\node at (-6,35) (h2) {$(2,1,3)$};
\node at (0,35) (h3) {$(2,3,0)$};
\node at (6,35) (h4) {$(3,0,2)$};
\node at (12,35) (h5) {$(3,1,0)$};

\node at (-6,40) (i1) {$(2,3,1)$};
\node at (0,40) (i2) {$(3,1,2)$};
\node at (6,40) (i3) {$(3,2,0)$};

\node at (0,44) (j) {$(3,2,1)$};


\draw[-, very thick,color=blue] (a) -- (b);

\draw[-] (b) -- (c1);
\draw[-, very thick,color=blue] (b) -- (c2);

\draw[-] (c1) -- (d1);
\draw[-] (c1) -- (d2);
\draw[-] (c1) -- (d3);
\draw[-] (c2) -- (d1);
\draw[-] (c2) -- (d3);
\draw[-, very thick,color=blue] (c2) -- (d4);

\draw[-] (d1) -- (e1);
\draw[-] (d1) -- (e2);
\draw[-] (d1) -- (e4);
\draw[-] (d2) -- (e1);
\draw[-] (d2) -- (e3);
\draw[-] (d3) -- (e2);
\draw[-] (d3) -- (e3);
\draw[-] (d3) -- (e5);
\draw[-, very thick,color=blue] (d4) -- (e4);
\draw[-] (d4) -- (e5);

\draw[-] (e1) -- (f1);
\draw[-] (e1) -- (f2);
\draw[-] (e1) -- (f3);
\draw[-] (e2) -- (f1);
\draw[-] (e2) -- (f2);
\draw[-] (e2) -- (f4);
\draw[-] (e2) -- (f5);
\draw[-] (e3) -- (f2);
\draw[-] (e3) -- (f6);
\draw[-] (e4) -- (f3);
\draw[-, very thick,color=blue] (e4) -- (f4);
\draw[-] (e4) -- (f5);
\draw[-] (e5) -- (f5);
\draw[-] (e5) -- (f6);

\draw[-] (f1) -- (g1);
\draw[-] (f1) -- (g2);
\draw[-] (f1) -- (g4);
\draw[-] (f2) -- (g2);
\draw[-] (f2) -- (g3);
\draw[-] (f2) -- (g6);
\draw[-] (f3) -- (g1);
\draw[-] (f3) -- (g3);
\draw[-] (f3) -- (g4);
\draw[-] (f3) -- (g6);
\draw[-, very thick,color=blue] (f4) -- (g1);
\draw[-] (f4) -- (g3);
\draw[-] (f4) -- (g5);
\draw[-] (f5) -- (g4);
\draw[-] (f5) -- (g5);
\draw[-] (f5) -- (g6);
\draw[-] (f6) -- (g6);

\draw[-] (g1) -- (h1);
\draw[-, very thick,color=blue] (g1) -- (h2);
\draw[-] (g2) -- (h1);
\draw[-] (g2) -- (h3);
\draw[-] (g2) -- (h4);
\draw[-] (g3) -- (h1);
\draw[-] (g3) -- (h3);
\draw[-] (g3) -- (h5);
\draw[-] (g4) -- (h2);
\draw[-] (g4) -- (h3);
\draw[-] (g4) -- (h4);
\draw[-] (g5) -- (h2);
\draw[-] (g5) -- (h3);
\draw[-] (g5) -- (h5);
\draw[-] (g6) -- (h4);
\draw[-] (g6) -- (h5);

\draw[-] (h1) -- (i1);
\draw[-] (h1) -- (i2);
\draw[-, very thick,color=blue] (h2) -- (i1);
\draw[-] (h2) -- (i2);
\draw[-] (h3) -- (i1);
\draw[-] (h3) -- (i3);
\draw[-] (h4) -- (i2);
\draw[-] (h4) -- (i3);
\draw[-] (h5) -- (i2);
\draw[-] (h5) -- (i3);

\draw[-, very thick,color=blue] (i1) -- (j);
\draw[-] (i2) -- (j);
\draw[-] (i3) -- (j);

\begin{scope}[red]
\node at (0,4) {\footnotesize{$(0,1)$}};
\node at (-1.5, 7.5) {\footnotesize{$(1,2)$}};
\node at (1.5, 7.5) {\footnotesize{$(0,1)$}};
\node at (3, 13.5) {\footnotesize{$(1,2)$}};
\node at (7, 13.5) {\footnotesize{$(0,1)$}};
\node at (1, 13.5) {\footnotesize{$(0,2)$}};
\node at (-5.6, 11.5) {\footnotesize{$(0,1)$}};
\node at (-5.5, 13.5) {\footnotesize{$(0,2)$}};
\node at (-3, 13.5) {\footnotesize{$(2,3)$}};
\node at (-10, 16.5) {\footnotesize{$(2,3)$}};
\node at (-9.5, 18.5) {\footnotesize{$(0,1)$}};
\node at (-7, 18.5) {\footnotesize{$(1,2)$}};
\node at (-3.5, 18.5) {\footnotesize{$(0,1)$}};
\node at (-1, 18.5) {\footnotesize{$(0,3)$}};
\node at (2,16.5) {\footnotesize{$(2,3)$}};
\node at (3.5, 19) {\footnotesize{$(0,1)$}};
\node at (6.5, 18.5) {\footnotesize{$(0,2)$}};
\node at (9, 18.5) {\footnotesize{$(0,2)$}};
\node at (11, 18.5) {\footnotesize{$(1,2)$}};
\node at (-14, 23) {\footnotesize{$(1,2)$}};
\node at (-12, 23) {\footnotesize{$(1,3)$}};
\node at (-10.5, 23.5) {\footnotesize{$(1,3)$}};
\node at (-8.5, 23.5) {\footnotesize{$(2,3)$}};
\node at (-7, 24) {\footnotesize{$(0,1)$}};
\node at (-4, 24) {\footnotesize{$(0,1)$}};
\node at (-1.5, 24) {\footnotesize{$(2,3)$}};
\node at (1, 24) {\footnotesize{$(0,1)$}};
\node at (3.5, 24) {\footnotesize{$(0,2)$}};
\node at (6, 24) {\footnotesize{$(1,2)$}};
\node at (8, 24) {\footnotesize{$(1,2)$}};
\node at (10, 24) {\footnotesize{$(0,1)$}};
\node at (12, 24) {\footnotesize{$(0,3)$}};
\node at (14, 23.5) {\footnotesize{$(2,3)$}};
\node at (-15, 28) {\footnotesize{$(0,1)$}};
\node at (-13, 28.5) {\footnotesize{$(0,2)$}};
\node at (-12, 29.5) {\footnotesize{$(0,3)$}};
\node at (-12.5, 27) {\footnotesize{$(2,3)$}};
\node at (-9, 29) {\footnotesize{$(1,2)$}};
\node at (-5, 29) {\footnotesize{$(0,1)$}};
\node at (-3, 29) {\footnotesize{$(0,3)$}};
\node at (-1, 28) {\footnotesize{$(2,3)$}};
\node at (.5, 29.5) {\footnotesize{$(0,2)$}};
\node at (1.5, 28.5) {\footnotesize{$(1,2)$}};
\node at (4, 29) {\footnotesize{$(1,3)$}};
\node at (7, 29) {\footnotesize{$(1,2)$}};
\node at (9, 27) {\footnotesize{$(0,1)$}};
\node at (11.5, 27) {\footnotesize{$(2,3)$}};
\node at (11.5, 28.5) {\footnotesize{$(1,3)$}};
\node at (11.5, 29.5) {\footnotesize{$(0,3)$}};
\node at (15, 28) {\footnotesize{$(0,1)$}};
\node at (-14, 33) {\footnotesize{$(2,3)$}};
\node at (-11.5, 33.5) {\footnotesize{$(0,1)$}};
\node at (-10, 34) {\footnotesize{$(0,2)$}};
\node at (-8, 34) {\footnotesize{$(1,2)$}};
\node at (-5, 34) {\footnotesize{$(0,1)$}};
\node at (-3, 34.5) {\footnotesize{$(0,3)$}};
\node at (-5, 32.5) {\footnotesize{$(0,2)$}};
\node at (-2, 31.5) {\footnotesize{$(1,2)$}};
\node at (.5, 34) {\footnotesize{$(0,3)$}};
\node at (3.5, 33) {\footnotesize{$(1,3)$}};
\node at (5.5, 34) {\footnotesize{$(2,3)$}};
\node at (7.5, 34) {\footnotesize{$(1,2)$}};
\node at (9.5, 34) {\footnotesize{$(1,3)$}};
\node at (11.5, 34) {\footnotesize{$(2,3)$}};
\node at (13.5, 32.5) {\footnotesize{$(0,1)$}};
\node at (-9, 38) {\footnotesize{$(1,2)$}};
\node at (-6, 39) {\footnotesize{$(1,3)$}};
\node at (-1, 36) {\footnotesize{$(0,1)$}};
\node at (1, 36) {\footnotesize{$(2,3)$}};
\node at (-2.5, 39.5) {\footnotesize{$(1,3)$}};
\node at (-1, 38.5) {\footnotesize{$(2,3)$}};
\node at (1, 38.5) {\footnotesize{$(0,1)$}};
\node at (2.5, 39.5) {\footnotesize{$(0,2)$}};
\node at (6, 39) {\footnotesize{$(0,2)$}};
\node at (9, 38) {\footnotesize{$(1,2)$}};
\node at (-3, 42) {\footnotesize{$(2,3)$}};
\node at (0, 42) {\footnotesize{$(1,2)$}};
\node at (3, 42) {\footnotesize{$(0,1)$}};
\end{scope}

\end{tikzpicture}
}
\caption{EL-labeling  of the rook monoid $R_3$.}
\label{F:Rook3}
\end{center}
\end{figure}

\begin{Theorem}\label{T:MainTheorem}
Let $\varGamma=\{0,1,\dots,n\} \times \{0,1,\dots,n\}$ 
and let $F:C(R_n)\longrightarrow \varGamma$ be the 
edge-labeling as defined as in (\ref{E:rooklabel}). 
Then $F$ is an EL-labeling  for $R_n$. 
\end{Theorem}
We will prove this theorem in the next section. 
The complete edge-labeling of $R_3$ is shown in Figure \ref{F:Rook3}.

\section{\textbf{Proofs}}\label{S:proofs}


Let $[x,y]$ be an interval in $R_n$, and let $\mathfrak{c}: x=x_0\lessdot \cdots \lessdot x_k=y$ be a maximal chain in $[x,y]$. 
The sequence $F(\mathfrak{c}) := (F(x_0,x_1), \dots, F(x_{k-1},x_k))$ will be called the 
{\em Jordan-H\"{o}lder sequence of $\mathfrak{c}$}. Here, $F$ is as defined in (\ref{E:rooklabel}).

\begin{Example}
Let $x=(0,1,0)$ and $y=(3,1,2)$ be two elements from $R_3$. 
It is easy to check from Figure~\ref{F:Rook3} that the maximal chain 
\[
\mathfrak{c}:\ x \lessdot (1,0,0) \lessdot (1,0,2) \lessdot (1,2,0) \lessdot (1,2,3) \lessdot (2,1,3) \lessdot y
\]
is the lexicographically smallest Jordan-H\"{o}lder sequence in $[x,y]$. 
Note that 
\[
F(\mathfrak{c}) = ((0,1),(0,2),(0,2),(0,3),(1,2),(2,3))
\]
is a non-decreasing sequence.  
\end{Example}

We will describe the lexicographically first chain, denoted by $\mathfrak{c}: x=x_0\lessdot \cdots \lessdot x_k=y$,    
in $[x,y]$. Clearly, if $F(x,x_1) = (a,b)$, then \hbox{$a=\min \{ c:\ (c,d)=F(x,z)\text{ for $z\in [x,y]$ with $x\lessdot z$}\}$.}
First, we assume that the least element in $x$ such that $a_i \neq b_i$ is nonzero.

\begin{Lemma}\label{L:short1}
If $a_i$ is the least element in $x$ 
such that $a_i \notin\{0, b_i\}$, then $x_1$ is obtained from $x$ by one of the following 
covering moves: 
\begin{itemize}
\item a type 2 covering relation where 
$a_i$ is interchanged with $a_i + s+ 1$ if $a_i +s+1$ is after $a_i$ but $a_i+1,\cdots, a_i+s$
appear in $x$ before $a_i$; 
\item a type 1 covering relation that is performed at $a_i$. 
\end{itemize} 
\end{Lemma}

\begin{proof}
We start with a claim.
Let $x=(a_1,a_2,\dots, a_n)$ and $y=(b_1,b_2,\dots, b_n)$ be the one-line representations for $x$ and $y$.
If $a_i$ denotes the least element such that $a_i \notin\{0, b_i\}$, then $a=a_i$.

First of all, it is easy to check by using the main result of~\cite{CR11} 
(that is the Deodhar-description of the Bruhat-Chevalley-Renner ordering on $R_n$) that $a_i < b_i$.
Secondly, if $a_i > a$, then $a_i\neq 0$, and therefore, $x$ and $y$ differ from each other at 
their nonzero entries only. But this immediately yields a contradiction with the minimality of $a_i$. 

We now assume towards a contradiction that $a_i< a$. 
If the number $b_i$ appears in $x$ as an entry $a_j$, then we have $i<j$. 
In this case, we have $a_j =b_i< b_j$ as well. 
We proceed with this assumption. 
Then the interchange of $a_i$ and $a_j$ gives us an element $z$ in $[x,y]$. 
But this means that there is another element $z'$ such that $x\lessdot z' \leq z$. 
Since $z$ agrees with $x$ on its first $i-1$ entries, so does $z'$. Furthermore, $z'$ 
can be chosen so that it is obtained from $x$ by moving its $i$-th entry.  
Then $F(x,z')$ is of the form $(a_i,b')$ for some $b'$. This contradicts the minimality of $(a,b)$. 

Let us now assume that $b_i$ does not appear in $x$. 
In this case, we define $z':=(c_1,\dots, c_n)$ by setting $c_j$ as $a_j$ for $j\neq i$ 
and by setting $c_i$ as the smallest number that is less than or equal to $b_i$ and not in $x$.
In this case, $z'$ covers $x$ with $F(x,z')=(a_i,c_i)$. Since $(a_i,c_i) < (a,b)$, 
we found another contradiction. This finishes the proof of our claim. 
Notice that the constructions of $z'$ in the two cases show also that 
the covering $x\lessdot x_1$ with lexicographically first label is obtained from $x$ as we stated.
This finishes the proof of our lemma. 
\end{proof}

We now consider the situation where the least element in $x$ such that $a_i \neq b_i$ is 0.
We have some related notation:
\begin{align*}
P(x,y):= \{ i:\ a_i=0, 0\notin \{a_{i+1},b_i\}\}\ \text{ and } \ 
Q(x,y):=\{ a_{i+1}:\ a_i=0, 0\notin \{a_{i+1},b_i\}\}.
\end{align*}
We set $q:=\min Q(x,y)$ and $s:= \min (\{1,\dots, n\} \setminus \{a_1,\dots, a_n\})$. 
Finally, let $a_{r+1}$ be the element in $Q(x,y)$ such that $q=a_{r+1}$. 

\begin{Lemma}\label{L:short2}
If $q\leq s$, then $x_1$ is obtained from $x$ by interchanging $a_r$ with $a_{r+1}$.
If $q > s$, then we have the following subcases.
\begin{enumerate}
\item If the number of zeros in $x$ and $y$ are the same, then 
$x_1$ is obtained from $x$ by interchanging $a_r$ with $a_{r+1}$.
\item If the number of zeros in $x$ and $y$ are different, then we replace $a_p$ with $s$,
where $p= \max P(x,y)$. 
\end{enumerate}
\end{Lemma}
\begin{proof}
The fact that these constructions give an element $x_1$ such that $x \lessdot x_1 \leq y$ 
is easily checked from Lemma~\ref{L:PPRcovering0}. 
The minimality of the resulting label, that is $F(x,x_1)$, follows from the construction of $x_1$. 
\end{proof}

Lemmas~\ref{L:short1} and \ref{L:short2} show that 
the lexicographically first chain in $[x,y]$ is the one that climbs up 
in the Bruhat-Chevalley-Renner order by moving the smallest ``out of place'' element in $x_i$ (for $i\in \{1,\dots, k\}$) 
to the left or by suitably incrementing it.
It is also clear from the construction of $x_1$, hence of $x_2,x_3$, and so on, that 
for each $i$ in $\{1,\dots, k-1\}$ there is only one covering relation whose label is the lexicographically smallest 
and such that the corresponding label is greater than or equal to the label that is constructed 
for the previous covering relation. In other words, 
there is a unique chain $\mathfrak{c}: x \lessdot x_1 \lessdot x_2\lessdot \cdots \lessdot x_k=y$ such that 
$F(x,x_1)\leq F(x_1,x_2) \leq F(x_2,x_3)\leq \cdots \leq F(x_{k-1},y)$. 
We summarize this observation as follows.

\begin{Proposition}\label{P:lexicographically}
Let $[x,y]$ be an interval in $R_n$. Then there is a unique maximal chain 
$\mathfrak{c}: x=x_0< \cdots < x_k=y$ in $[x,y]$ 
whose Jordan-H\"{o}lder sequence $F(\mathfrak{c})$ is the lexicographically first sequence. 
Furthermore, this chain is increasing in the sense that 
\begin{equation*}
F(x_0,x_1) \leq F(x_1,x_2) \leq \cdots \leq F(x_{k-1},x_k).
\end{equation*}
\end{Proposition}

\begin{proof}[Proof of Theorem \ref{T:MainTheorem}.]
Let $F:C(R_n)\longrightarrow \varGamma$ denote the edge-labeling that is defined in (\ref{E:rooklabel}). 
By Proposition~\ref{P:lexicographically}, we know that there is a unique lexicographically 
first maximal chain in $R_n$. In other words, $F$ is an EL-labeling, and hence, $R_n$ is an EL-shellable poset. 
\end{proof}

It is well-known that an interval of length two in a Weyl group consists of four elements. 
However, for the Renner monoids this is not true in general. In $R_n$, already for $n=2$, 
there are intervals of length two with three elements. 
We finish this section by proving that in $R_n$ an interval of length two has at most four elements.

\begin{Proposition}\label{C:chainordiamond}
Let $[x,y] \subseteq R_n$ be an interval of length two. 
Then, $[x,y]$ is either a chain, or a diamond. In other words, either  
$[x,y]=\{x,x_1,x_1',y\}$ with $x< x_1 \neq x_1' < y$, or $[x,y]=\{x,x_1,y\}$ with $x< x_1  < y$.
\end{Proposition}

\begin{proof}
Let $(a_1,\dots,a_n)$ and $(b_1,\dots,b_n)$ denote $x$ and $y$, respectively. 
Let $J$ denote the set $\{ i \in [n]:\ a_i \neq b_i \}$. 
If $\ell(y) - \ell(x) =2$, then by Lemmas \ref{L:PPRcovering0} and \ref{L:PPRcovering1}, we see that $|J| \leq 4$. 
Therefore, if $z=(c_1,\dots, c_n) \in [x,y]$ is strictly between $x$ and $y$, then the set 
$\{ i \in [n]:\ c_i \neq a_i \ \text{or}\ c_i \neq b_i \}$ has at most four elements as well.
Arguing case by case, as in Lemmas \ref{L:short1} and \ref{L:short2},
we see that there are at most two different possibilities for $z$. 
\end{proof}

\section{\textbf{Final remarks}}\label{S:Final}

Let $P$ be a poset, and let $\widehat{P}$ denote 
the poset that is obtained from $P$ by adjoining to it 
$\hat{0}$ (the smallest element) and $\hat{1}$ (the maximal element). 
We denote by $I(P)$ the 
set of all intervals in $\widehat{P}$. 
The M\"{o}bius function $\mu: I(P) \longrightarrow \Z$ is an integer valued function, 
(uniquely) determined by the following conditions: 
\begin{equation*}
\mu ([ x,y ]) =
\begin{cases}
1 & \, \text{if $x=y$}; \\
-\sum_{x \leq z < y} \mu([ x,z]) & \, \text{if $x < y$}.
\end{cases}
\end{equation*}
Let $ R_{n,k} \subset R_n$ ($0\leq k \leq n$) denote the subposet consisting of elements whose rank is $k$. 
In \cite{ACT}, it is shown that the M\"obius function on $I(R_{n,k})$ takes values in $\{-1,0,1\}$. 
When $k=n$, $R_{n,k}$ is the symmetric group, and the M\"obius function on $S_n$ is well-known
(see~\cite{Verma71, Stembridge07, Jones09}).

The {\em order complex} of $P$, denoted by $\Delta(P)$, is the 
abstract simplicial complex whose faces are the chains in $P$. 
There is a remarkable topological interpretation 
of the M\"obius function of $\widehat{P}$ in relation with $\Delta(P)$. 
For $x,y\in P$, we have $\mu ((x,y)) = \widetilde{\chi} ( \text{lk}_{x,y} )$, 
where $\text{lk}_{x,y}$ is the order complex of the open interval 
$(x,y) = \{ z \in P:\ x < z < y \}$, 
and $\widetilde{\chi}$ is the reduced Euler characteristic of the topological space $\text{lk}_{x,y}$.
In particular, $\mu ( \widehat{P} ):= \mu ([\hat{0},\hat{1}])$ is the reduced Euler characteristic of $\Delta(P)$. 
In our next result we will determine the homeomorphism type of $\Delta(R_n)$. To this end, let us recall 
the following important result of Danaraj and Klee.
\begin{Lemma}[\cite{DK74}, page 444]\label{L:DK}
Let $\Delta$ denote a pure shellable simplicial complex in which every codimension one face is contained in at most two facets.
If every codimension one face is contained in exactly two facets, then $\Delta$ is homeomorphic to a sphere. Otherwise,
$\Delta$ is homeomorphic to a ball. 
\end{Lemma}

\begin{Theorem}\label{T:final}
The order complex $\Delta ( (x,y) )$ of every open interval $(x,y) \subset R_n$ is homeomorphic to a sphere, or a ball. 
In particular, $\Delta(R_n)$ triangulates a ball of dimension $n^2$.
\end{Theorem}
\begin{proof}
It follows from Proposition \ref{C:chainordiamond} that every codimension one face of an open interval $(x,y)$ in  $R_n$
is contained in at most two facets. Thus, by Theorem \ref{T:MainTheorem} and Lemma \ref{L:DK} we see that $\Delta ( (x,y) )$
has the homotopy type of a sphere or a ball.

For our second claim, first we observe that the rank of $\hat{1}$ in $R_n$ is $n^2$ by Corollary~\ref{C:length for Sn}.
Hence, $\dim \Delta(R_n) = n^2$.
To see that $\Delta(R_n)$ is homeomorphic to a ball, it suffices to notice that $\Delta(R_n)$ has a unique atom,
so it is a cone.
\end{proof}

\begin{Remark}
Recently, a remarkable generalization of Theorem~\ref{T:final} is obtained by Warut Thawinrak 
in his Undergraduate Honors Thesis at the University of Minnessota-Twin Cities. 
He describes an effective and easy method for deciding whether or not the 
M\"obius function on an interval $[x,y]\subseteq R_n$ takes 
only nonzero values, and in the non-vanishing situation, he shows that the interval $[x,y]$ is 
 isomorphic to an interval in a suitable symmetric group.
\end{Remark}

\bibliography{Shellable}
\bibliographystyle{plain}
\end{document}